\documentclass[final]{siamltex}
\usepackage{graphicx,amssymb}

\newcommand{\de}{\partial}
\newcommand{\De}{\Delta}
\newcommand{\cF}{\mathcal{F}}

\newcommand{\lu}[1][l]{\;^#1\!}
\def\norm#1{\Vert #1 \Vert_1}
\def\normi#1{\Vert #1 \Vert_\infty}
\newcommand{\w}{w_j^{(i)}\,}
\newcommand{\et}{\mathcal{E}}
\newcommand{\wu}{\bar{w}^{(i)}\;}
\newcommand{\ldu}{\lu\bar{\lambda}}

\newcommand{\vu}{\bar{v}^{(n)}}
\newcommand{\Ep}{\hat{E}}
\newcommand{\Ed}[1][j]{\;\lu E^i_{#1}}
\newcommand{\Rp}{\hat{R}}
\newcommand{\cO}{\mathcal{O}}

\title{A finite difference method for Piecewise Deterministic Processes with memory}

\author{Mario Annunziato
\thanks{Dipartimento di Matematica ed Informatica, Universit\`a degli Studi di
Salerno, Via Ponte Don Melillo, 84084 Fisciano (SA), Italy 
({\tt mannunzi@unisa.it}). {\sl Preprint n. 11-2006}.}
}

\begin{document}
\maketitle
\begin{abstract}
In this paper the numerical approximation of solutions of Liouville-Master 
Equations for time-dependent distribution functions of Piecewise Deterministic Processes 
with memory is considered. 
These equations are linear hyperbolic PDEs with non-constant
coefficients, and boundary conditions that depend on integrals over the 
interior of the integration domain.
We construct a finite difference method of the first order, by a combination
of the upwind  method, for PDEs, and by a direct quadrature, 
for the boundary condition.
We analyse convergence of the numerical solution for distribution functions
evolving towards an equilibrium.
Numerical results for two problems, whose analytical solutions
are known in closed form, illustrate the  theoretical finding.
\end{abstract}

\begin{keywords}
Piecewise-deterministic process, dichotomic noise,
random telegraph process, binary noise, upwind method,
conservative systems.
\end{keywords}

\begin{AMS}
65-06, 65M06, 65M12, 60K40, 60J25, 60J75
\end{AMS}

\section{Introduction}
We deal with the following system of PDEs:
\begin{equation}
      \de_t F_s(x,y,t) + A_s(x)\, \de_x F_s(x,y,t) + \de_y F_s(x,y,t) =
      - \lambda_s(y)\,F_s(x,y,t)
\label{basic_kolmogorov}
\end{equation}
with Cauchy initial conditions:
\begin{equation}
   F_s(x,y,t_0) = F_{0,s}(x)\delta(y)
%
   \label{basic_cauchy}
\end{equation} 
and boundary conditions:
\begin{equation}
  F_s(x,0,t) = \sum_{j=1}^{S} q_{sj} \int_{0}^{t-t_0}\!dy\; 
  F_j(x,y,t)\,\lambda_j(y)
  \label{basic_boundary}
\end{equation}
for the $s=\left\lbrace 1, \ldots, S \right\rbrace$ unknowns
 $F_s:\mathcal{D}\rightarrow\mathbb{R}$,
with $(x,y,t) \in \mathcal{D} := (\Omega\times [0,T-t_0]\times [t_0,T]) \subset \mathbb{R}^3$, where 
$\Omega = [\Omega_a,\Omega_b] \subset \mathbb{R}$. 
$q_{sj}$ are the elements of a stochastic matrix 
having the following fundamental properties:
$0 \leq q_{sj}\leq 1$ and $\sum_s q_{sj} =1$.
The known functions $F_{0,s}(x)$, $A_s(x)$ and $\lambda_s(y)\geq 0$,
will be discussed later. 

Eq. (\ref{basic_kolmogorov}), jointly with boundary conditions, 
is the Liouville-Master equation 
\footnote{In general, equations for density probability of random processes
are derived from the Chapman-Kolmogorov equation. As discussed in Ref. \cite{gar}
the same equation turns into a Liouville equation in absence of randomness, and
into a Master Equation, if only jump processes are involved. We use both terms
in order to stress the deterministic and the random character of the processes
considered here.} 
for the probability distribution functions $F_s(x,y,t)$ of
a continuous \textit{piecewise-deterministic process} (PDP), that has
been introduced by Davis \cite{dav_bk,dav}.
Indeed, here we deal with a simplified version of Davis' PDPs, but still enough 
general to cover many interesting models. 
The definition of PDP is more popular between researchers 
working on operations research and probability calculus 
(see, e.g., \cite{cos:duf}), rather than others outside these fields, 
even though the latter unknowingly use it, at least in a simplified form. 
Before to proceed with the discussion of the 
numerical solution of our problem, we give a short introduction of the 
underling PDP process we are considering here.
\footnote{The author acknowledges Prof. M.H.A. Davis for some explanations about
the definition of PDP.}
\begin{definition}
\label{def:PDP}
We name $X(t)$, $X:\mathbb{R}\rightarrow\mathbb{R}$, be a continuous PDP if:
\begin{itemize}
\item[(a)] 
$X(t)$ satisfies the equation:
\begin{equation}
	\dot{X}(t)=A_s(X),   \;\;\;\;\;\;\;\; s=1,\ldots,S
\label{RTE_estesa_1}
\end{equation}
where $A_{s}:\mathbb{R}\rightarrow\mathbb{R}$ is a function chosen randomly 
on a set of $\lbrace A_1,\ldots,A_S\rbrace$ known functions. 
Given $A_s$, we say that the dynamics is in the (deterministic) state $s$. 
We require that $A_s(x)$ be Lipschitz 
continuous, so that, for fixed $s$, $X(t)$ exists, is unique and non-explosive
solution.
\item[(b)]
The initial condition is settled by the 
Cauchy problem to Eq. (\ref{RTE_estesa_1}), 
i.e. $X(t_0)=X_0$, and by the initial state $s=s_0$ of the same equation. 
\item[(c)] 
Each state $s$ is characterised by an its own probability density function (PDF)
$\psi_s:\mathbb{R}^+\rightarrow\mathbb{R}^+$ of ``transition events'': 
\begin{equation}
\psi_s(t), \;\mbox{ with }  \int_0^\infty \!\! \psi_s(t) dt =1.
\label{RTE_estesa_2}
\end{equation}
\item[(d)] 
When a transition event occurs, the dynamics switches instantaneously from the 
state $j$ ($A_j$) to a new state $i$ ($A_i$), given randomly by the
transition probability matrix (or transition measure):
\begin{equation}
   q_{ij}
\label{RTE_estesa_3}
\end{equation}
The position $X(t)$ of the process is not affected when the state switches.
\end{itemize}
\end{definition}
Assumptions (\ref{RTE_estesa_1}), (\ref{RTE_estesa_2}) and (\ref{RTE_estesa_3})
define the three \textit{local characteristics} of the PDP.
We see that Eq. (\ref{RTE_estesa_1}) can be 
integrated as an ordinary differential equation, provided that not any switching 
event happen inside the integration interval. Therefore, with the exceptions of 
switching times, the process is deterministic, continuous and composed 
of pieces of solutions of Eq. (\ref{RTE_estesa_1}).
Anyway, the whole resulting process $X(t)$ is not
deterministic, it represents a random sample path on a probability
space.
\footnote{For our pourposes, we do not need to specify what the abstract 
probability space is.}

The statistical description of the process is given by the
unknown functions $F_s(x,y,t)$: each represents the probability
to find the process $X(t)$, in the state $s$, at time $t$ in a position less than $x$,
being past the time $Y$ since the last switching event.
Formally we write:
$$
F_s(x,y,t) := \mathbb{P}(X(t)\leq x, y \leq Y < y+dy, \mbox{state}=s)
$$
where $\mathbb{P}$ is a probability measure of a  
probability space for the process.
If we are interested only in the position $x$ of the process,
we can integrate over all values of $y$, and the distribution function 
for the process, regardless the time $y$, reads as:
\begin{equation}
	\mathcal{F}_s(x,t) := \int_0^{t-t_0}\!\! dy\;F_s(x,y,t).
\label{dist_senza_memoria}
\end{equation}
With the further hypothesis $X(t)\in \Omega$, we have:
\begin{equation}
\de_x F_s(x,y,t)\vert_{x=\de\Omega} =0
\label{neumann}
\end{equation}
since there is a null probability for the process to be outside
the interval $\Omega$.
Besides, the probability measure have to be conserved during the evolution, so that:
\begin{equation}
  \lim_{x \rightarrow \Omega_b} \sum_{s=1}^S \mathcal{F}_s(x,t) = 1,
\;\;\;\;\;\;\;\forall t
\label{conservation+}
\end{equation}
and
\begin{equation}
  \lim_{x \rightarrow \Omega_a} \sum_{s=1}^S \mathcal{F}_s(x,t) = 0,
\;\;\;\;\;\;\;\forall t
\label{conservation-}
\end{equation}
have to be satisfied. This three last equations are boundary conditions
for (\ref{basic_kolmogorov}), that complete the definition of the
problem we approach to treat here.

The function $\lambda(y)$, named \textit{hazard function} (or \textit{hazard rate}), 
is related to the statistics of the PDF switching times (\ref{RTE_estesa_2}) by:
\begin{equation}
  \lambda_s(y):=\frac{\psi_s(y)}{\int_y^\infty \psi_s(\tau)\,d\tau}.
\label{hazard}
\end{equation}
It represents the probability per unit of time that a transition
event will occur, i.e. a transition rate, having past the time $y$ 
since the last event.
The explicit dependence of $\lambda$ on $y$ makes both the statistics
of the switching events and the process $X(t)$ be non-Markovian, 
so that $y$ plays the role of \textit{memory}.

The main aim of this article is to solve Eq. (\ref{basic_kolmogorov}), jointly
to all the above mentioned boundary conditions, by a finite difference
scheme of the first order and prove convergence of the numerical solution. 
We note that numerical methods for solving 
linear hyperbolic PDEs with non-constant coefficient, such as (\ref{basic_kolmogorov}), 
are well known in literature \cite{mor:may,ran}, but what makes
this problem a little special is the boundary condition (\ref{basic_boundary}):
the value of the unknowns $F_s$, on the boundary $y=0$, depend on
an integration of the $F_s$ over the interior of the domain.
This means that the numerical scheme for Eq. (\ref{basic_kolmogorov})
have to be supported by one for (\ref{basic_boundary}), so that
conditions for convergence of the numerical solution 
have to be investigated again. As a result we found a Courant-Friedrichs-Lewy 
(CFL) condition for ensuring linear convergence.

The secondary, but not of minor importance, aim of this article is to provide 
a connection bridge between PDPs as known by experts of the field and, as above mentioned,
the same processes as known by others, who apply them to modelling in several
areas of science and engineering.
Here we give a sample of quotas, for which PDPs can be concerned by others,
grouped in two categories: diffusive processes and systems having an equilibrium. 
We mention:
anomalous diffusion \cite{bol:gri:wes}, reaction-diffusion \cite{hor}, 
scattering of radiation  \cite{jak:ren}, biological dispersal \cite{oth:dun:alt}, 
for the former category,
and non-Maxwellian equilibriums \cite{mario,grigo,cac,fil:hon,kit,bro:han},
diagnostic techniques for semiconductor lasers \cite{jak:rid}, 
filtered telegraph signals \cite{paw:ric,jak:rid},  harmonic
oscillators \cite{mas:por}, ecological systems \cite{man:etal},  for the latter.
Many of the models involved in such references, concern the application
of a two-state noise to a dynamical equation.
\footnote{For more related citations, the reader can search the following 
key words: dichotomic/binary noise/process, random telegraph process, coloured noise.} 

The common end of all these researches, consists in extracting statistical 
properties from process governed by that equation. 
Generally an approximation method can be applied to the original model:
such as by the projector technique \cite{zwa,mori.h,grigo},  by
a ``coarse-grain'' technique,  (see, e.g., \cite{hil:oth}), an asymptotic analysis
(see, e.g., \cite{bol:gri:wes}).
However not always these techniques provide a satisfactory description.
In some cases an exact analytical result can obtained as in Refs. \cite{paw:ric,jak:rid,mor}
and more recently by the characteristic functional method \cite{bud:cac}.
Obviously, computations can also be performed by Monte Carlo's simulations, 
but, at the best of our knowledge, few or nothing it has been devoted to 
a direct calculation of the time-dependent distribution function comprensive of an 
explicit memory variable.
Concerning this, we remark that the main alternative 
is based on the inclusion of \textit{supplementary variables} \cite{dav_bk,cox2}, 
that turns PDP into Markovian, i.e. a memoryless process.

In the next section we provide an example that emphatizes the connection 
between PDPs and models with dichotomic noise, and a conjecture that ensures the
existence of a stationary solution of Eq. (\ref{basic_kolmogorov}).
In Sect. \ref{sec:num_scheme} we establish the numerical scheme.
We introduce definitions in Sect. \ref{sec:conv} and in Sect. \ref{sec:analysis}
some theoretical results about the related convergence.
In Sect. \ref{sec:comp} we present  numerical results to two problems 
for which an analytical stationary solution is known in closed form, and
verify the stated convergence properties.

\section{Explanatory example}
\label{sec:rc-filter}

Let us consider a dissipative process $X(t)$ subject to a noised input $\xi(t)$, described
by the equation:
\begin{equation}
\frac{dX}{dt} = - X(t) +  \xi(t).
\label{FRTP}
\end{equation}
If $\xi(t)$ is taken as the random telegraph signal, Eq. (\ref{FRTP}) act as filter, 
and $X(t)$ is referred as filtered random telegraph process \cite{paw:ric,jak:ren}.
The same equation is elsewhere referred as Langevin equation \cite{mario,cac} 
subject to a dichotomous noise.
$\xi(t)$ alternately takes on values $\pm 1$, 
with an exponential (or Poisson) statistics for the transition events (\ref{RTE_estesa_2}): 
 $\psi(\tau) = \mu \, e^{-\mu \tau}$, where $\mu^{-1}$ is the
average time $\langle\tau\rangle$ between transitions.
The process $X(t)$ results composed of pieces of increasing and decreasing exponentials.
The statistical properties of the process $X(t)$ can be found
by the associated probability density distributions $p^{\pm}(x,t)$,
governed by a Liouville-Master Equation \cite{won,jak:rid,mario}:
\begin{equation}
\left\{
\begin{array}{l}
\de_t p^+ -(x-1) \de_x p^+ = (1-\mu) p^+ + \mu p^-  \\
\de_t p^- -(x+1) \de_x p^- = \mu p^+ + (1-\mu) p^-.
\end{array}
\right.
\label{ME_Filter_RTS}
\end{equation}
Now let us to see the same process from the point of view of PDPs.
The  exponential statistics for $\psi(t)$ makes the process of transitons be Markovian and 
the hazard function constant: $\lambda(t) = \mu$. 
Eq. (\ref{basic_kolmogorov}) turns into:
$$
      \de_t F_s(x,y,t) + A_s(x)\, \de_x F_s(x,y,t) + \de_y F_s(x,y,t) =
      - \mu \,F_s(x,y,t).
$$
By integrating this equation over all the values of $y$, we get:
$$
      \de_t \mathcal{F}_s(x,t) + A_s(x)\, \de_x \mathcal{F}_s(x,t) - F_s(x,0,t) 
      = - \mu \,\mathcal{F}_s(x,t),
$$
having used the property $F_s(x,y,t)=0$ if $y>0$, since the process
is memoryless. 
From Eq.(\ref{basic_boundary}) we have:
$$
 F_s(x,0,t) = \mu \sum_{j=1}^{S} q_{sj}  \int_{0}^{t}\!dy\; F_j(x,y,t)= \\
 		\mu \sum_{j=1}^{S} q_{sj}  \cF_j(x,t)
$$ 
and inserting it into the previous:
$$
      \de_t \cF_s(x,t) + A_s(x)\, \de_x \cF_s(x,t)  =
      \mu \sum_{j=1}^{S} (q_{sj}-\delta_{sj}) \cF_j(x,t)
$$
If $S=2$, with the known functions:
\begin{equation}
A_1(x)=(1-x), \;\;\;\;\;\;\;\;\; A_2(x)=-(1+x),
\label{field_example}
\end{equation}
and transition measure:
\begin{equation}
q_{11}=q_{22}=0, \;\;\;\;\;\;\;\; q_{12}=q_{21}=1,
\label{trans_matrix_example}
\end{equation}
provided that $p^{\pm}(x,t)=\de_x \cF_{1,2}(x,t)$,
we obtain just the equation (\ref{ME_Filter_RTS}). 
This show the connection between PDPs and processes driven by dichotomous noise.

\subsection{Remarks on equilibrium solutions}

In what follows we focus our attention on solutions $F(x,y,t)$
having an equilibrium, but we presume that the numerical scheme can be extended 
to diffusion processes too.
Conditions for the existence of equilibrium solution can be conjectured
by using simple dynamical arguments \cite{mor,ben:etal}.
If all dynamical equations (\ref{RTE_estesa_1}) own only attraction points and
all these are contained into the intersection of the basin of 
attraction of each, then a process starting from this region
will never escape. Whence, there should exists a region $\Omega$ where the 
process is confined and a stationary distribution 
$\cF_{eq}(x)=\lim_{t\rightarrow\infty}\cF(x,t)$ exists.

\section{The finite-difference scheme}
\label{sec:num_scheme}

In this section we show the numerical scheme to solve Eqs. (\ref{basic_kolmogorov})
and (\ref{basic_boundary}) based on a finite difference method of first order.
For the shake of simplicity we take $t_0=0$ and the domain of $F_s$ becomes 
$\mathcal{D}:= \{\Omega\times [0,T]\times [0,T]\}$.
It is convenient to perform the numerical integration along 
the characteristic lines $\xi = t-y$. With this new variable, we define the unknowns 
 $ \phi_s(x,y,\xi)= \phi_s(x,y,t-y) := F_s(x,y,t) $, so that 
Eq. (\ref{basic_kolmogorov}) transforms as:
\begin{equation}
	A_s(x)\, \partial_x\phi_s(x,y,\xi) + \partial_y \phi_s(x,y,\xi) =
	-\lambda_s(y) \phi_s(x,y,\xi).
\label{basic_kolm_charact}
\end{equation} 
This equation is valid for $0 < \xi < t$ and $0 < y < t$. The initial condition is given on 
\begin{equation}
	\phi_s(x,y,\xi)\vert_{\xi=-y} = F_{0,s}(x)\delta(y)
\label{basic_cauchy_charact}
\end{equation}
With the new variable we get $\phi_s(x,0,\xi)|_{\xi=t} = F_s(x,0,t)$, and
the boundary condition Eq. (\ref{basic_boundary}) becomes: 
\begin{equation}
\phi_s(x,0,t) = \sum_{l=1}^{S} q_{sl} \int_0^{t} 
  dy \; \phi(x,y,\xi)|_{\xi=t-y} \lambda_l(y)
\label{basic_boundary_charact}
\end{equation} 
We will assume that similar conditions of Eqs. (\ref{dist_senza_memoria}), 
(\ref{conservation+}) and (\ref{conservation-}) are satisfied for
$\phi_s(x,y,\xi)$, and a stationary solution exists.

\newcommand{\fid}[1][kj]{\lu\phi_{#1}^{i}}
\newcommand{\ld}{\lu\lambda_j}
\newcommand{\Ad}{\lu A_k}
On the domain $\mathcal{D}$, we introduce a uniform mesh:
\begin{equation}
 (x_k,y_j,t_n) \;\;\;\;
\left\{
\begin{array}{l}
	k=0,\ldots, N_k \\
	j,n=0,\ldots,N, \;\;\;\;\; N=T/\De t
\end{array}
\right.
\label{uniform_mesh}
\end{equation}
with step size $\De x$ and $\De y =\De t$, so that we define the discrete 
known functions as  $ \lu A_{k} := A_l(x_k)$ and $\ld := \lambda_l(y_j)$,
and the discrete solution:
$$
\lu F_{kj}^n,  \;\;\;\;\;\; n=0,\ldots, N,  \;\;\;\; j<n
$$
as an approximation of $F_l(x_k,y_j,t_n)$ at the mesh points.

The change of variable $\xi = t-y$ corresponds to the following discrete mapping
on the mesh:
\begin{equation}
(k,j,n)=(k,j,i)|_{i=n-j}, 
\label{mapping_discreto}
\end{equation}
therefore we get the following relation:
\begin{equation}
F_l(x_k,y_j,t_n) = \phi_l(x_k,y_j,t_n-y_j)=\phi_l(x_k,y_j,\xi_i) \Rightarrow
						\lu F_{kj}^n = \lu \phi_{kj}^{n-j} = \fid
\label{map_tra_le_soluzioni}
\end{equation}
between the discrete solutions.
Here the index $i$ identifies the set of mesh points lying on
the characteristics lines.

The numerical scheme is obtained by discretizating both 
equations (\ref{basic_kolm_charact}) and (\ref{basic_boundary_charact}).
We apply upwind to the first, and get:
\begin{eqnarray}
\label{basic_kolm_disc}
\fid[k,j+1] = \fid - \Ad \frac{\De y}{\De x}
      \left(\fid[k+\nu,j] - \fid[k+\nu-1,j] \right) - \ld\; \fid \De y  \\
\;\;\;\; i=0,\ldots,N   \nonumber
\end{eqnarray} 
where $\nu=1$ if $\Ad < 0$, and $\nu=0$ if $\Ad > 0$.
The boundary condition (\ref{neumann}) is included by requiring
that $\fid[0j]=\fid[1j]$ and $\fid[N_k-1,j]=\fid[N_k,j]$.

For the second, we substitute 
integral with a quadrature scheme:
\newcommand{\fidi}{\lu\phi_{k\,j}^{i-j}\;}
\begin{equation}
\lu[s]\phi_{k,0}^i\; = \De y \sum_{l=1}^S q_{sl} \sum_{j=0}^{i} \w\;\; \fidi\; \ld 
\;\;\;\;\;\;\; i>0
\label{basic_boundary_disc} 
\end{equation} 
where $\w \geq 0$ is a sequence of weights.

\begin{figure}
\begin{center}
\includegraphics[scale=0.8]{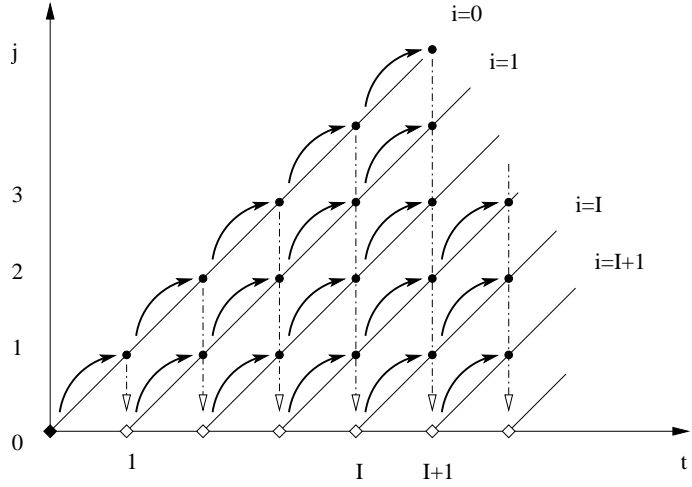}
\end{center}
\caption{Representation of the integration scheme on the regular mesh of the
$(t,y)$ plane. Full curved arrows: upwind step of Eq. (\ref{basic_kolm_disc}).
Dot-dashed vertical arrows: quadrature of Eq. (\ref{basic_boundary_disc})}
\label{fig_integ_scheme}
\end{figure}

The integration proceeds as follows. Given the initial condition 
$\lu\phi_{k,0}^0=\phi_l(x_k,0,0)\\=F_l(x_k,0,0)$, Eq. (\ref{basic_kolm_disc}) 
allow us to calculate all the
$\fid$ starting from $j=1$ up to $N$, for the fixed characteristic line $i=0$.
In general given the values on the boundary $j=0$: $\fid[k,0]$,
we can find all $\fid$ starting from $j=1$ up to $N-i$, for a fixed characteristic line $i$
(see curved arrows of Fig. (\ref{fig_integ_scheme})).
But starting values $\fid$ for upwind are unknown
and have to be estimated by using the boundary integration
(\ref{basic_boundary_disc}) (see vertical dot-dashed arrows of Fig. (\ref{fig_integ_scheme})).
This is a system of equations for the unknowns $\fid[k,0]$.
When all $\fid$ are known, the discrete distribution function can be retrieved by 
\newcommand{\Fd}{\lu F_{k,j}^{n}\;}
\newcommand{\FFd}{\lu\cF_{k}^{n}\;}
$ \Fd = \fid|_{i=n-j}$, and also the discrete countrepart
of (\ref{dist_senza_memoria}):
\begin{equation}
\FFd:=\sum_{j=0}^{n} v^{(n)}_j  \Fd \Delta y,
\label{dist_no_y_disc}
\end{equation} 
can be estimated by a quadrature formula of weights $v^{(n)}_j$.

%
\newcommand{\ed}[1][k,j]{\;\lu e_{#1}^{i}}
\newcommand{\Ki}{K}

\section{Preliminary definitions}
\label{sec:conv}

\subsection{Global errors and convergence}

We are interested in how well $\fid$, $\Fd$ and $\FFd$ approximate the corresponding
analytical solutions.
We consider the global pointwise error:
\begin{equation}
   \ed := \phi_{l}(x_k,y_j,\xi_i) - \fid
\label{global_err}
\end{equation}
for the transformed solution. The same value defines the global error for 
$\Fd$ under the discrete mapping (\ref{mapping_discreto}).

In order to prove convergence, we introduce norms for measuring errors. 
For spatial $x_k$ and memory variable  $y_j$ the $\infty-$norm is used.
The discrete $1$-norm for the states $s$ of the system is the natural choice, because
of the conservation of the probability of Eq. (\ref{conservation+}).
For convenience of notation we define the global error for the state
$l$ at time $t_i$ and time memory $y_j$ as:
\begin{equation}
  \Ed := \max_k|\ed|
\end{equation}
and the global error regardless states as:
\begin{equation}
\norm{E_j^i} := \sum_{l=1}^S \Ed
\label{err-stato}
\end{equation}

We say that $\Fd$ converges to $F_l(x_k,y_j,t_n)$ in the norm $\Vert\cdot\Vert$ if:
\begin{equation}
  \Vert E^i\Vert = \max_j \norm{E_j^i} \rightarrow 0, \;\;\;\; 
	\mbox{as } \De x,\De y \rightarrow 0.
\end{equation}

The global error for the distribution function (\ref{dist_no_y_disc}) is defined as:
\begin{equation}
   \Ep_{k}^n := \cF(x_k,t_n) - \cF_k^n = 
	\int_0^{t_n} dy \sum_l F_l(x_k,y,t_n) -  
	\De y \sum_{j=0}^n v_j^{(n)} \sum_l \Fd
\label{distrib_err}
\end{equation}
and the associated convergence is stated by:
\begin{equation}
  \normi{\Ep^n} = \max_k |\Ep^n_k| \rightarrow 0, \;\;\;\;
	\mbox{as}\; \De x,\De y \rightarrow 0.
\label{err_distr}
\end{equation}

\subsection{Local truncation error and quadrature error}

As usual \cite{mor:may,ran}, the local truncation error is defined 
by inserting the true solution $\phi_l(x,y,\xi)$ into the 
discrete scheme of Eq. (\ref{basic_kolm_disc}), i.e.:
\begin{eqnarray}
\lu\varepsilon_{kj}^i  := && \frac{\phi_l(x_k,y_i+\De y,\xi_i)-\phi_l(x_k,y_i,\xi_i)}{\De y} \nonumber \\
\label{loc_trunc} && - A_l(x_k)\frac{\phi_l(x_k+\nu \De x,y_i,\xi_i) - 
						\phi_l(x_k+(\nu-1) \De x,y_i,\xi_i)}{\De x} \\ 
&& - \lambda_l(y_j) \phi_l(x_k,y_i,\xi_i). \nonumber
\end{eqnarray}
By evaluating the rest with respect to $\fid$, we get:
\begin{equation}
\lu\varepsilon_{kj}^i = \frac{1}{2} \De x \left(
|A_l(x_k)| \de_x^2 \phi_l(\eta_k,y_j,\xi_i) -
\alpha \de_y^2\phi_l(x_k,\eta_j,\xi_i)    \right)
\label{trunc_err}
\end{equation}
where $\alpha:=\De y/\De x$, and $\eta_k,\eta_j$ are unknown points.

The quadrature error committed from (\ref{basic_boundary_disc}) for 
the evaluation of the boundary integral (\ref{basic_boundary_charact}) is 
defined as:
\begin{equation}
\lu R_k^i := \mathcal{M}_i t_i \De y^\beta \de_y^\beta 
(\phi_l(x_k,\tilde{\eta}_j ,\xi_i)\, \lambda_l(\tilde{\eta}_j))
\label{quad_err}
\end{equation}
where $\mathcal{M}_i$ are some constants that can depend on $i$, and
$\tilde{\eta}_j$ are unknown points  of the local integration interval.
$\beta$ defines the order of repeated quadrature 
formulas (\ref{basic_boundary_disc}).

The quadrature error committed from (\ref{dist_no_y_disc}) for 
the evaluation of (\ref{dist_senza_memoria}) is defined as:
\begin{equation}
	\lu \Rp_k^n := \hat{\mathcal{M}}_n t_n \De t^{\hat{\beta}} 
	\de_y^{\hat{\beta}} F_l(x_k,\hat{\eta}_j,t_n),
\label{quad_err_no_y}
\end{equation}
where, as for the previous error, $\hat{\mathcal{M}}_n$ are some values that can depend on $n$, and
$\hat{\eta}_j$ are unknown points.
$\hat{\beta}$ defines the order of (\ref{dist_no_y_disc}).

\section{Analysis of convergence}
\label{sec:analysis}

In this section we show first a lemma for convergence of the 
numerical solution $\fid$ for the transformed equation (\ref{basic_kolm_disc}),
then proof a theorem for the convergence order of the numerical
solution $\cF_k^n$ to the distribution function $\cF(x_k,t_n)$.
Proofs are based on classical arguments by finding bounds for global errors.

\begin{lemma}
\label{lm:convergenza}
Let  $\phi_l(x,y,\xi) \in C^{2,\bar{\beta},2}(\mathcal{D})$ be a solution 
of  Eq. (\ref{basic_kolm_charact}) under the boundary conditions (\ref{basic_cauchy_charact}), with $\bar{\beta}=\max\{\beta,2\}$, $\max_l\normi{A_l(x)}\leq M$,
and \\ $\max_l\normi{\lambda_l(y)} \leq L_u$, 
for $x\in \Omega$ and $l\in \{1,\ldots, S\}$.
Let  $(x_k,y_j,t_n)$ be an uniform mesh on $\mathcal{D}$ defined in 
(\ref{uniform_mesh}), of step sizes $\De x$ and $\De y =\De t$. 
Let $\fid$ be the numerical solution resulting from the scheme as defined in 
Eqs. (\ref{basic_kolm_disc}) and (\ref{basic_boundary_disc}), under the
transformed mesh of Eq. (\ref{mapping_discreto}).

If the Courant-Friedrichs-Lewy (CFL) condition 
\begin{equation}
\De y < \left(\frac{M}{\De x}+L_u\right)^{-1}
\label{CFL}
\end{equation}
is satisfied then:
\begin{enumerate}
\item[1.] Given the error $\norm{E_0^i}$ at boundary $y=0$, the error $\norm{E_m^i}$ computed at time step $t_{i+m}$ along the characteristic $\xi_i$, is bounded
by:
\begin{equation}
\norm{E_m^i} \leq \norm{E_0^i} (1-L_u \De y)^m + m\De y \norm{\et^i}
\label{tesi_2}
\end{equation}
where $\norm{\et^i}=\cO(\De x), \;\;\forall i$.
\item[2.]  If $\De y^{-1}>\max_i(w_0^{(i)}) \max_l \lambda_l(0)$, 
there exist constants $L_d$, $\Ki$, 
such that the error computed at time $t_i$ along the boundary
$y=0$ is bounded by:
\begin{equation}
	\norm{E^i_0} \leq \norm{\bar{R}} \exp\left(\Ki t_i\right)
 	+ \norm{\bar{\et}} \Ki \frac{t_{i+1}^2}{2}
	\exp \left( (\Ki - L_d) t_i \right)
\label{tesi_1}
\end{equation}
where  $\norm{\bar{R}}=\cO(\De x^\beta)$ and $\norm{\bar{\et}}:=\max_i \norm{\et^i}=\cO(\De x)$.
\end{enumerate}
\end{lemma}

This lemma states that the numerical solution $\fid$ converges to analytical 
with a linear order.
Errors calculated for $\phi$ are same for $F$.
Being interested in finding the probability regardless of
memory and state at fixed time $t_n$, we search an estimate for the  error
$\normi{\Ep^n}$ as defined in Eq. (\ref{err_distr}).
\begin{theorem}
\label{th:conv}
Let the hypothesis of the Lemma \ref{lm:convergenza} be satisfied,
then for the problem of Eq. (\ref{basic_kolmogorov})
the numerical solution (\ref{dist_no_y_disc}), obtained as discussed
in Sec. \ref{sec:num_scheme}, converges to analytical
(\ref{dist_senza_memoria}) in the sense (\ref{err_distr}), with order $\cO(\De x)$.
\end{theorem}

\begin{proof}
We substitute the integration of Eq. (\ref{dist_senza_memoria}) over the memory state
with a sequence of quadrature formulas of weights $v_j^{(n)}$:
\begin{equation}
\int_0^{t_n} dy \sum_l F_l(x_k,y,t_n) = 
\De y \sum_{j=0}^n v_j^{(n)}  \sum_l F_l(x_k,y_j,t_n)+ \sum_l \lu \Rp_k^n 
\label{quadrature_y}
\end{equation}
so that:
$$
\normi{\Ep^n} = \max_k \left| 
   \De y  \sum_j v_j^{(n)} \sum_l F_l(x_k,y_j,t_n)+ \sum_l \lu \Rp_k^n
         - \De y \sum_j v_j^{(n)}  \sum_l \Fd
  \right|
$$

Let $\vu := \max_j |v_j^{(n)}|$ and $\Rp^n:=\sum_l \max_k |\lu \Rp_k^n|$, 
from the triangle inequality:
$$
\normi{\Ep^n} \leq 
\vu \De y \sum_j \sum_l \max_k |F_l(x_k,y_j,t_n) 
- \Fd| + \Rp^n
$$
by using the discrete mapping (\ref{mapping_discreto}) and the definition
of the error (\ref{err-stato}), we get:
$$
\normi{\Ep^n} \leq \vu \De y\sum_{j=0}^n \norm{E^{n-j}_j}  + \Rp^n.
$$
By inserting Eq. (\ref{tesi_2}) we find:
$$
\normi{\Ep^{n}} \leq \vu \De y\sum_{j=0}^n 
	(\norm{E^{n-j}_{0}} (1-L_u \De y)^j + j \De y \norm{\et^i}) + \Rp^n
$$
and
$$
\normi{\Ep^{n}} \leq \vu \De y\sum_{m=0}^n 
\norm{E_0^m}(1-L_u\De y)^{n-m} + \vu \norm{\bar{\et}} \frac{t_{n+1}^2}{2} + \Rp^n
$$
where $\norm{\bar{\et}} := \max_j \norm{\et^{n-j}}$.
This inequality relates the searched probability error to errors along
the boundary $j=0$. Now we insert the second result of the previous lemma
stated by Eq. (\ref{tesi_1}) and find the order of the error for
vanishing $\De y$.
From the summation we get:
$$
\sum_{i=0}^n (1-L_u\De y)^{n-i} \left(
\norm{\bar{R}} e^{\Ki t_i}   
 	+ \norm{\bar{\et}} \Ki \frac{t_{i+1}^2}{2}
	e^{( \Ki - L_d) t_i}
\right)
$$
that is of order:
$$
\norm{\bar{R}}\frac{e^{K t_n}-1}{(K + L_u)\De y} +
\norm{\bar{\et}} K \frac{ t_n^2 e^{(K-L_d)t_n}}{2(K+L_u-L_d) \De y}.
$$
Finally, we get an order of convergence for the error:
\begin{eqnarray}
\label{error_order}
\normi{\Ep^{n}} \lesssim \;
   &&   \vu \norm{\bar{R}}\frac{\exp(K t_n)-1}{K+L_u} +  
      \Rp^n +  \\
   && \vu \norm{\bar{\et}} K \frac{ t_n^2 \exp((K-L_d)t_n)}{2(K+L_u-L_d)} +
      \vu \norm{\bar{\et}} \frac{t_n^2}{2}   \nonumber 
\end{eqnarray}
\end{proof}


\section{Computational results}
\label{sec:comp}

In what follows we present results of the numerical scheme, 
applied to the examples considered in Ref. \cite{paw}. 
For quadrature of Eqs. (\ref{basic_boundary_disc}) and (\ref{quadrature_y}),
we adopt the rectangle scheme:
\begin{equation}
\w = v_j^{(i)} = 
\left\{ 
	\begin{array}{l}
	0 \;\;\mbox{for } j=0  \\
	1 \;\;\mbox{for } j>0  \\
	\end{array}
\right.
\;\;\;\;\;\; i>0
\label{pesi}
\end{equation}
whose quadrature error (\ref{quad_err}) is:
$$
\lu R_k^i = \frac{1}{2} t_i \De y \de_y
(\phi_l(x_k,\tilde{\eta}_j ,\xi_i)\, \lambda_l(\tilde{\eta}_j )).
$$
Note that despite the low order ($\beta=1$) of approximation for quadrature,
the global error of the numerical solution is not degradated, because
the same order apply to upwind. The choice $w_0^{(i)}=0$ makes the equation
(\ref{basic_boundary_disc}) be explicit.
This can be consistent, because
for vanishing $\De y$ the contribution to integral is 
also vanishing for a limited $\fid[k,0]$.
For the explicit scheme the computational cost can be evaluated as follows:
at time step $i$ all upwinds take $2 N_k S i$ operations, the boundary 
quadrature takes $N_k S^2 i$. By summing over $i$ we get:
$$
 \mbox{Computational Cost} \approx N_k S^2 N^2.
$$

The Cauchy problem for starting the numerical integration is set, according to 
Eq. (\ref{basic_cauchy}), as follows:
\begin{equation}
\begin{array}{l}
\lu[s]F_{k0}^0 = 
	\left\{
	\begin{array}{lr}
		0, \;\;\;\; & k<0 \\
		0.5/(S \De y), \;\;\;\;& k=0 \\
		1  /(S \De y), \;\;\;\;&  k>0
	\end{array}
	\right. \\
 \\
\lu[s]F_{kj}^0  = 0, \;\;\;\; j > 0
\end{array}
\label{discrete_cauchy}
\end{equation}
for all $s = 1, \ldots, S$. This choice is an approximation of (\ref{basic_cauchy}) 
with $F_{0,s}(x)=H(x)$, where $H(x)$ is the Heaviside function.
Such Cauchy conditions for the Liouville-Master Equation correspond to having placed
the process $X(t)$ at the initial position $X_0=0$, to an equiprobable random initial state
, having spent the time zero (i.e. $\delta(y)$ in $(\ref{basic_cauchy})$).
This is a common choice when studying the motion of a particle
subject to a random fluctuating force, but is not a good mathematical hypothesis
for applying Lemma \ref{lm:convergenza}. However, it is well known the upwind 
methods tends to regularise the solution (numerical viscosity)
around discontinuities \cite{ran,mario:prep}, and, for the problems we are approaching
to solve, a unique stationary solution exists regardless the initial state of the process.
In this way we are enabled to use such ``non-regular'' Cauchy condition
for our numerical convergence tests. $\Fd$ are
then integrated by (\ref{quadrature_y}).

We are interested in plotting the density probability function of the 
process, regardless memory and states, defined as:
\begin{equation}
p(x,t) := \de_x \cF(x,t) := \de_x \sum_{l=1}^S \cF_l(x,t).
\label{pdf}
\end{equation}
The discrete version of this operation is:
\begin{equation}
p_k^n := \frac{1}{\De x} (\cF_{k+1}^n - \cF_{k}^n) := 
	\frac{1}{\De x} \sum_{l=1}^S (\lu \cF_{k+1}^n - \lu \cF_{k}^n)
\label{pdf_num}
\end{equation}
that is the first order right-derivative of the numerical distribution function.

\subsection{RC-filter subject to Markovian process (Poisson PDF)}
\label{sec:poisson}

\begin{figure}[!t]
\includegraphics[scale=0.5]{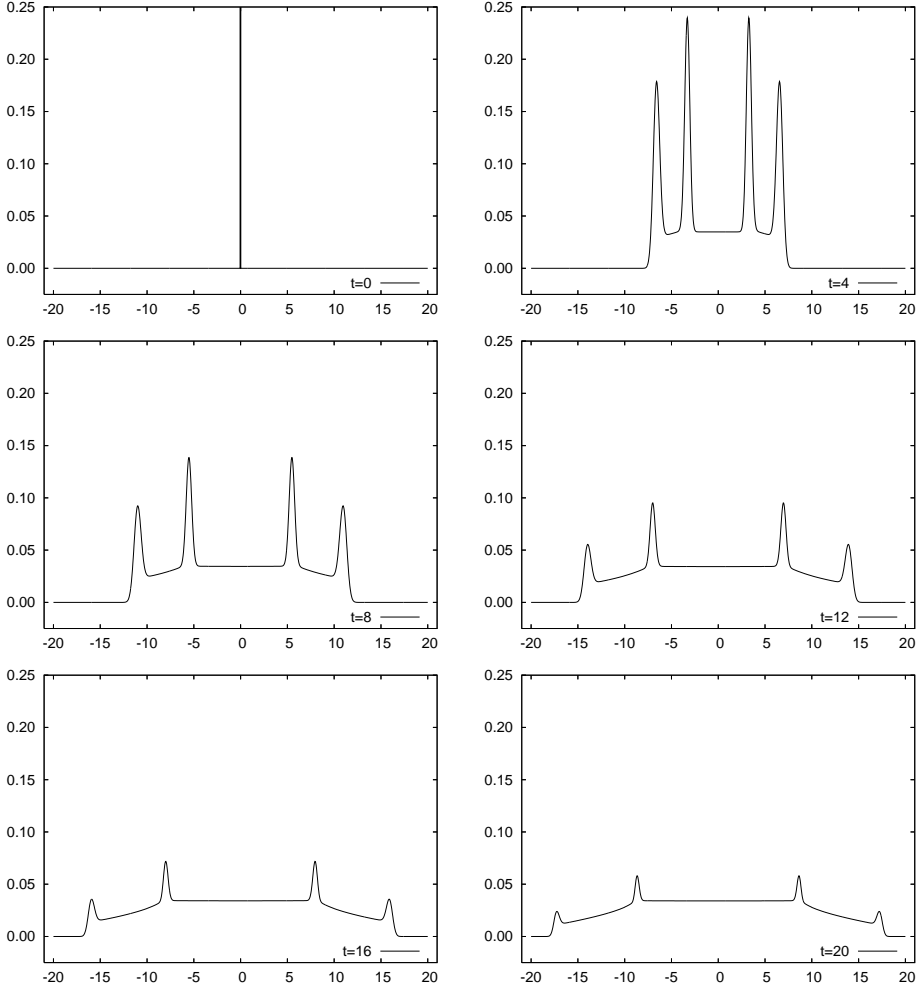}
\caption{Six snapshots of the temporal evolution of $p_k^n$ for the Langevin 
equation driven by  Poisson distribution time intervals of Sec. \ref{sec:poisson}
(horizontal axis: $x_k$; $t_n=0,4,8,12,16,20$).}
\label{fig:markov}
\end{figure}

In Fig. (\ref{fig:markov}) are plotted six snapshots of the 
temporal evolution related to the  four states process:
$A_s(x)=-\gamma_s x + W_s$, studied in \cite{mario:prep}.
Parameters are: $W_0=1,W_1=-1,W_2=2,W_3=-2$, $\lambda_s=0.2$ and 
$\gamma_s=0.1 \;\;\;\forall s$.
Results are comparable with that of the cited reference.
Note the broadening of the four peaks due to the numerical viscosity 
of the upwind method.

\subsection{Filtering of non-Markov dichotomous noise
with McFadden interval PDF}
\label{sec:mcfadden}

\begin{figure}[!t]
\includegraphics[scale=0.5]{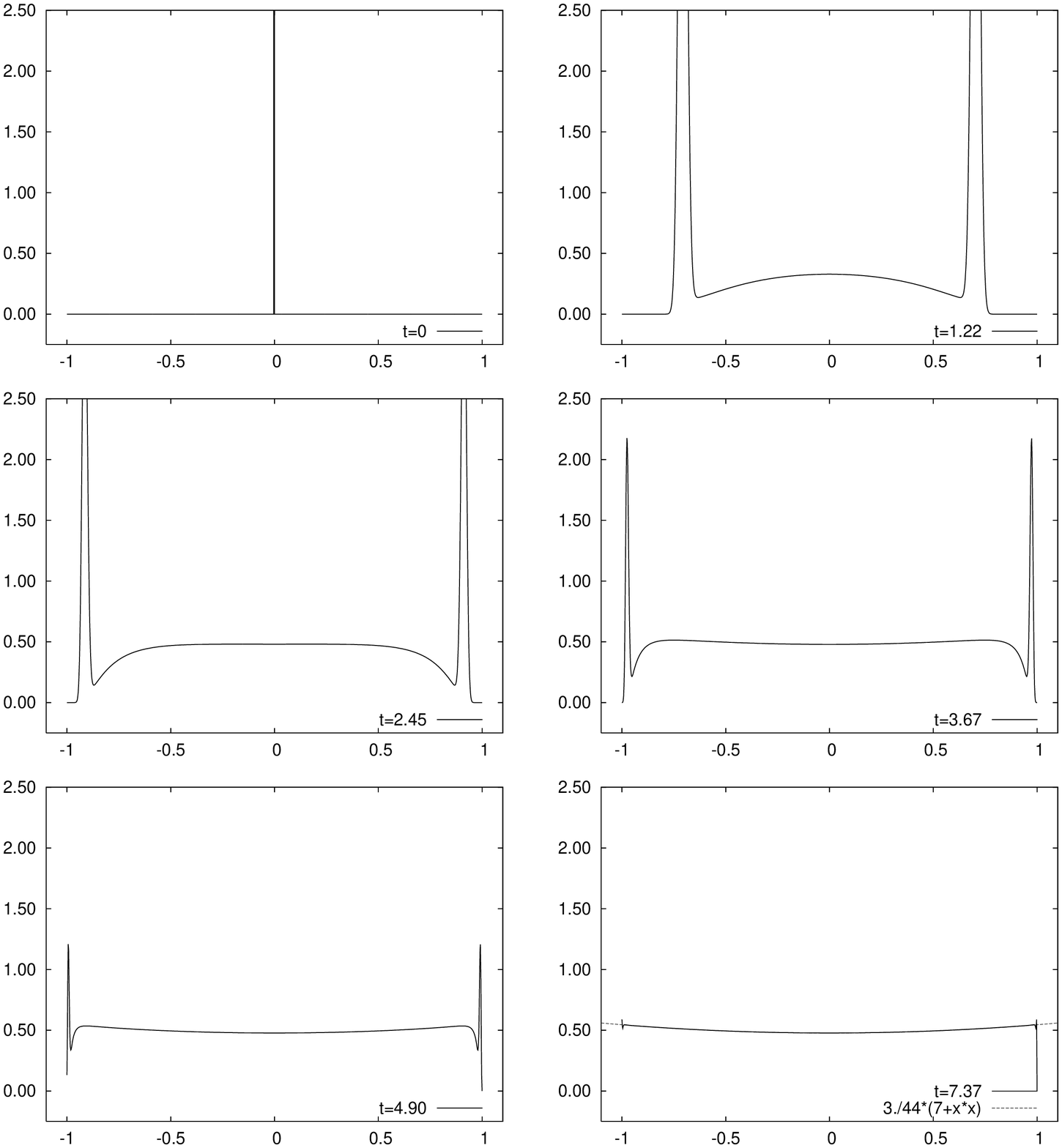}
\caption{Six snapshots of the temporal evolution of $p_k^n$  for 
filtering of dichotomous noise driven by the McFadden distribution time intervals
(horizontal axis: $x_k$).}
\label{fig:mcfadden}
\end{figure}

For this example the process is described by the Langevin
equation (\ref{FRTP}), with functions of Eq. (\ref{field_example})
and transition measure of Eq. (\ref{trans_matrix_example}).
Intervals between switching time have the McFadden PDF:
$
 \psi(t) = 3\,e^{-t}(1-e^{-t})^2
$.
The equilibrium density distribution has the form \cite{paw}:
\begin{equation}
p_{eq}(x)=\lim_{t\rightarrow\infty} p(x,t)= \frac{3}{44}(7+x^2),\;\;\;\; |x|<1
\label{mcfadden}
\end{equation}
and its integral:
\begin{equation}
\cF_{eq}(x) = \lim_{t\rightarrow\infty} \cF(x,t)= 
\frac{1}{44}(x^3+ 21 x +22),\;\;\;\; |x|\leq 1
\label{eq:soluz_mcfadden}
\end{equation}
The hazard function related to the density for 
switching intervals is:
$$
  \lambda(t)=\frac{3(1-e^{-t})^2}{2-e^{-t}+(1-e^{-t})^2}.
$$
We note that is $\lambda(0)=0$, so that the error committed from the choice  
$w_0^{(i)}=0$ (see (\ref{pesi})), is further improved.
Beside it is $\lambda(t) \leq 4/9$ and the convergence lemma give us
more guarantee that the errors does not grows fastly. 

Here the grid step sizes are $\De x= 0.002 $ and $\De y= 8\cdot 10^{-4}$.
Integration starts with a concentrated initial density (\ref{discrete_cauchy}) and 
stops at time $T=7.37$.

In Fig. \ref{fig:mcfadden} we see six snapshots of 
the numerical solution of the PDF (\ref{pdf_num}). At time $t=0$
the density of the process is concentrated to $x=0$, then two peaks,
corresponding to the two dynamical states, moves towards the
attraction points $x\pm 1$, and at last the stationary distribution appears.

\subsection{Filtering of non-Markov dichotomous noise
with ``gamma'' interval PDF}
\label{sec:gamma}

\begin{figure}[!t]
\includegraphics[scale=0.5]{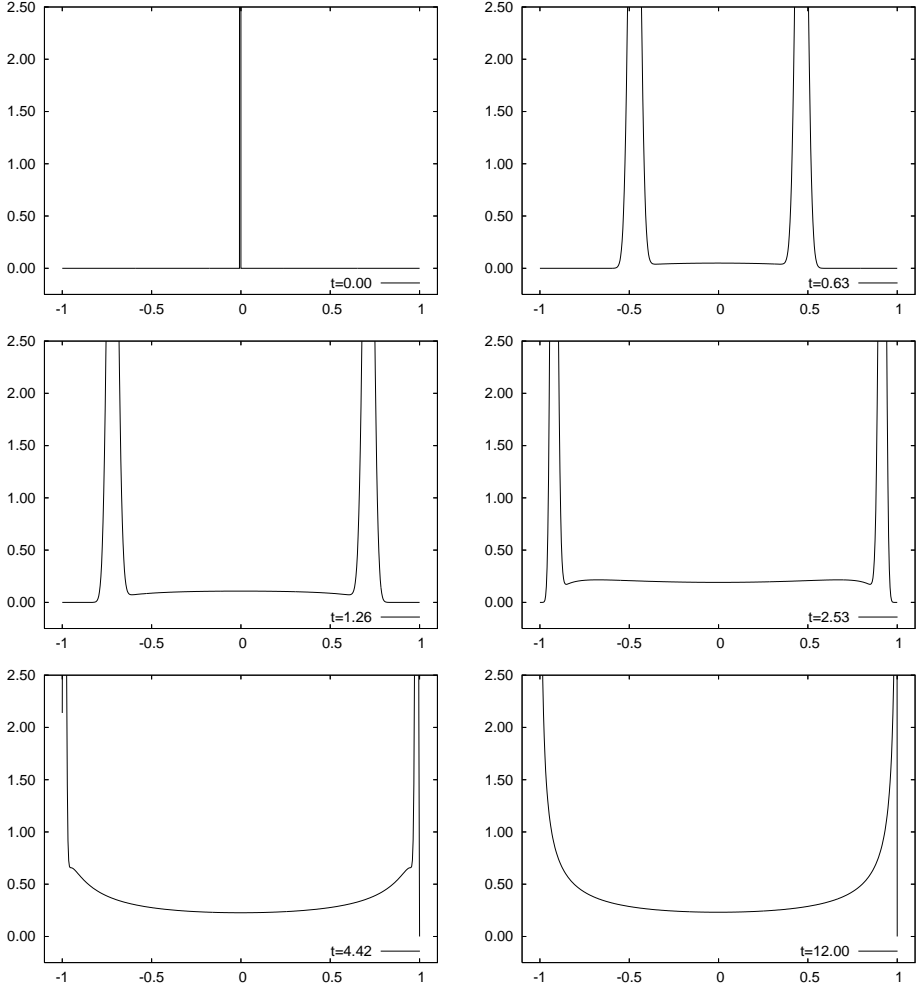}
\caption{Six snapshots of the temporal evolution of $p_k^n$  for 
filtering of dichotomous noise driven by the gamma distribution time intervals
(horizontal axis: $x_k$).}
\label{fig:gamma}
\end{figure}

For this example the process is described by the same 
Langevin equation as that the previous one,
but the intervals between switching
times of $\xi(t)$ of Eq. (\ref{FRTP}) are taken as the gamma density 
$
 \psi(t) = \mu^2t e^{-\mu t}
$
for both states \cite{paw}. Provided that $\mu=1/2$, the equilibrium
solution for the total density distribution function is:
\begin{eqnarray}
\label{eq:soluz_dist_gamma}
p_{eq}(x)=\lim_{t\rightarrow\infty} p(x,t) = 
\frac{|\Gamma(1-r)|^2}{\pi^{3/2}}(1-x^2)^{-1/2}\; 
	  _2F_1(r,r^*;\frac{1}{2};x^2),   \\
		r=(1+\sqrt{-1})/4, \;\;\; |x|\leq 1  \nonumber
\end{eqnarray}
in which $_2F_1$ is a hypergeometric function.

The \textit{hazard function} (\ref{hazard}) related to $\psi(t)$ is:
$$
   \lambda(t)= \frac{\mu^2 t}{1+\mu t}
$$ 
We see from the convergence theorem  that 
the error does not grow so fastly, because 
the maximum value of $\lambda(t)$ is $\lambda_{max}=\mu/e$. 
We perform the numerical integration on a mesh with spatial
discretization step $\De x = 0.004$ and temporal step $\De t = 0.0015$.
Integration starts with a  (\ref{discrete_cauchy}) and stops time 
$T=12$, where the equilibrium is 
supposed to be reached in good approximation.

In Fig. \ref{fig:gamma} are plotted six snapshots of the
total density distribution (\ref{pdf_num}). The evolution
behaves  as the previous example, with the exception of singularities
at $x\pm1$, for the equilibrium.


\begin{figure}[!t]
\includegraphics{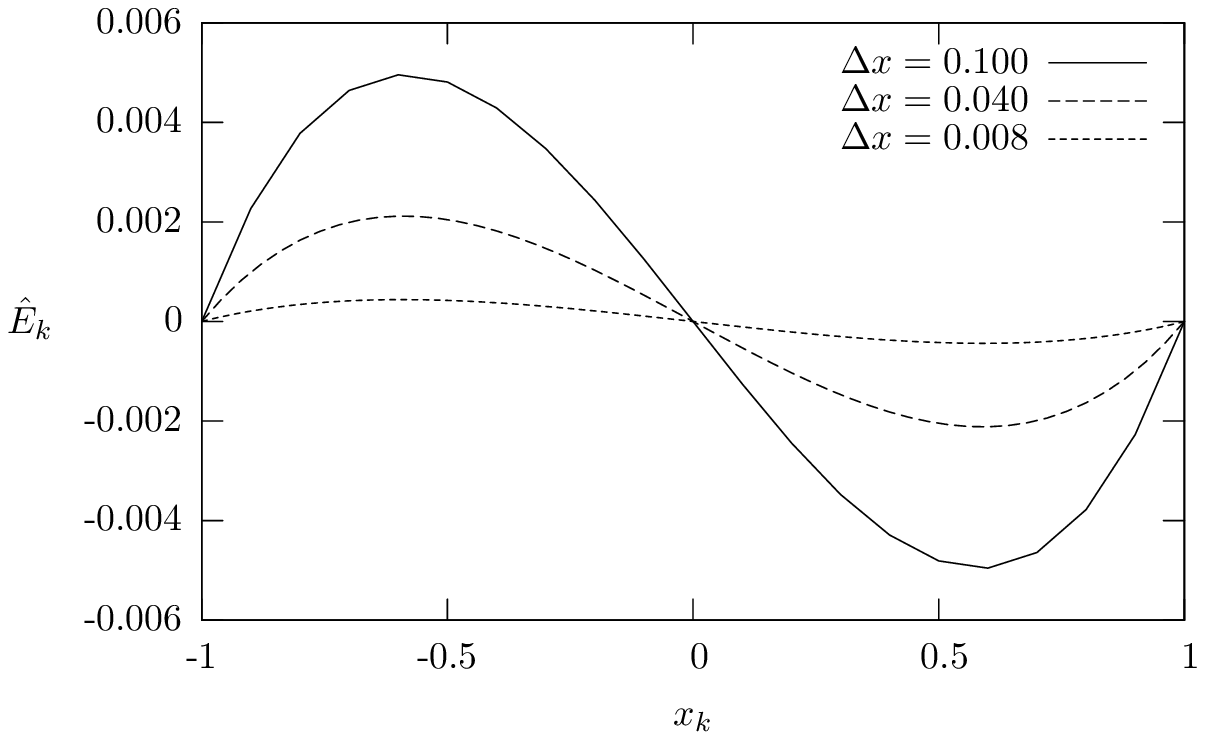}
\caption{Global error $\Ep_k$ between the exact solution of 
Eq. (\ref{eq:soluz_mcfadden}) (McFadden PDF) and 
numerical calculated for $\De x=0.1$,$\De x=0.04$,$\De x=0.008$.}
\label{fig:mcfad_errs}
\end{figure}

\begin{figure}[!t]
\includegraphics{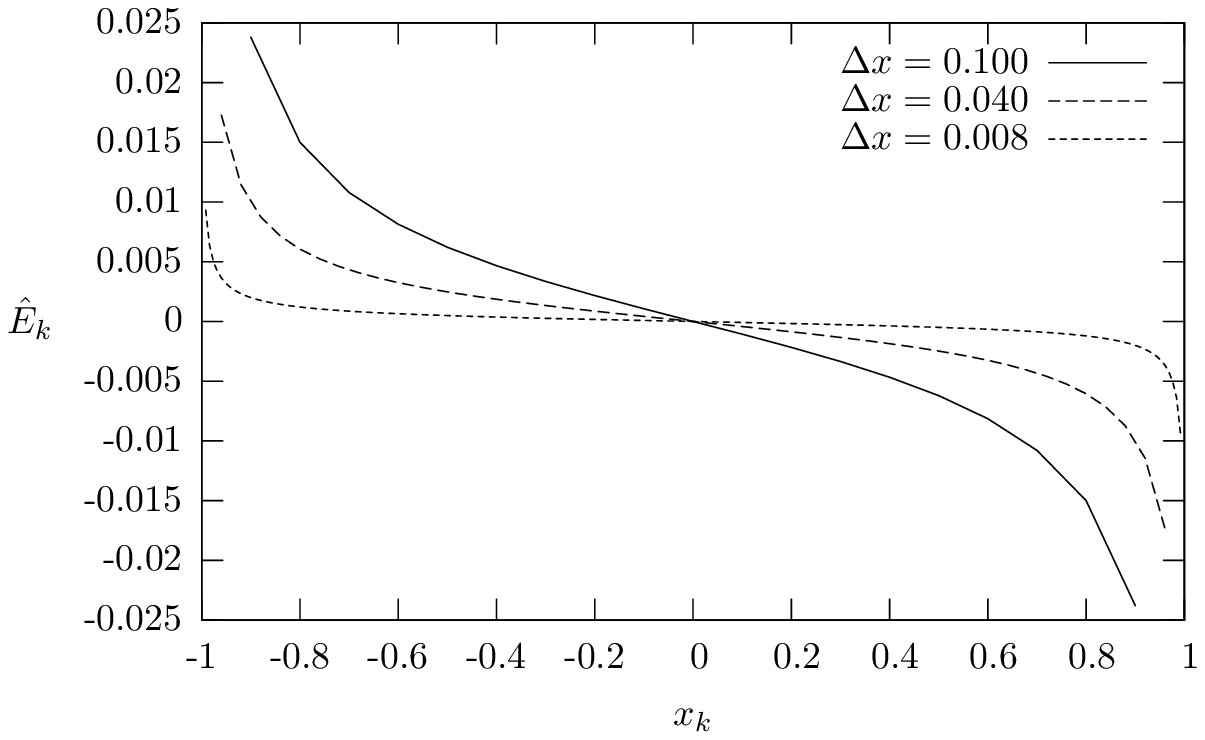}
\caption{Global error $\Ep_k$ between the integral of the exact solution of 
Eq. (\ref{eq:soluz_dist_gamma}) (gamma PDF) and 
numerical calculated for $\De x=0.1$, $\De x=0.04$, $\De x=0.008$.}
\label{fig:gamma_errs}
\end{figure}

\subsection{Convergence tests}

Since for the above mentioned problems we know two analytical results,
we can calculate the global error $\normi{\Ep}$ for the stationary
distribution functions of Eqs. (\ref{eq:soluz_mcfadden}) and (\ref{eq:soluz_dist_gamma}).
The analytical solution of Eq. (\ref{eq:soluz_dist_gamma}) is evaluated
by using MATLAB\circledR
\footnote{MATLAB\circledR \hspace{2pt} is a trade mark of ``The Matworks, Inc.''}
with libraries for calculating the Hypergeometric function \cite{hyper}.
We calculate the solution with the numerical scheme until the time $T=20$.
At this time we experienced that the stationary solution is reached.
This integration is repeated for some spatial step size $\De x$,
with the temporal step size constraint $\alpha=\De y/\De x = 0.9\,(M+L_u \De x)^{-1}$,
satisfying the CFL condition (\ref{CFL}).
In Fig. (\ref{fig:mcfad_errs}) we show the global error $\Ep_k$ plotted
for $\De x= \{0.1,0.04,0.008\}$ for the McFadden intervals. 
We see clearly that the maximum error decreases as $\De x$ decreases.
In Fig. (\ref{fig:gamma_errs}) we show the same test for gamma distributed intervals.
Also here the error decreases, but it shows a sort of divergence near $x=\pm 1$.

\begin{figure}[!t]
\includegraphics[scale=1.2]{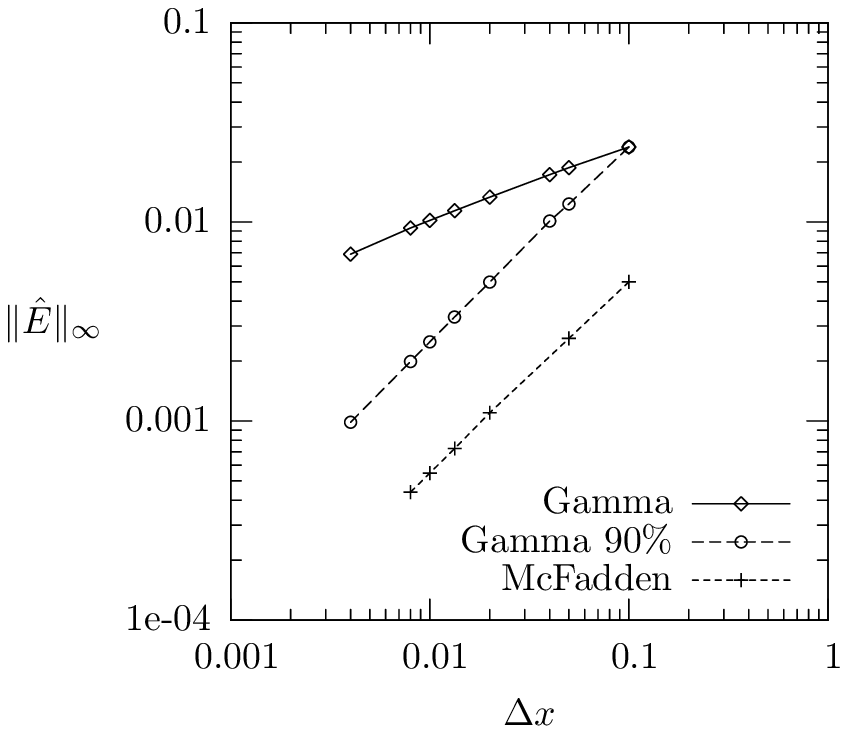}
\caption{Test for global error $\normi{\Ep}$ vs. mesh step size $\De x$.}
\label{fig:converge}
\end{figure}

In order to stress convergence, in Fig. (\ref{fig:converge}), we plot 
the error $\normi{\Ep}$ versus the step size $\De x$.
We see that the McFadden's  data's have unitary slope, i.e. linear convergence, 
as we expected from Theorem (\ref{th:conv}). Instead
gamma's data's are arranged with approximately $1/2$ slope.
This can be explained as follows.
We know \cite{paw} that the PDF of Eq. (\ref{eq:soluz_dist_gamma}) is $\cup$-shaped 
having $(1-x^2)^{-1/2} \ln(1-x^2)$ infinities near $x\pm 1$.
As $\delta x := 1-x$ approaches to zero, the second spatial 
derivative of $\cF_{eq}$ behaves as:
$$
\de_x p_{eq} = \de_x^2 \cF_{eq} \propto \delta x^{-3/2}(1 + \ln(2\,\delta x) ).
$$
Being $A_1(\delta x)\approx \delta x$ and by considering
$\De x = \delta x$, we get for the error:
$$
\normi{\Ep} \approx \sqrt{\De x} (1 + \ln(2\,\De x) ),
$$
that explain what we see from the results.
If we remove the endpoints from the measure of the error, 
we see from ``gamma 90\%'' of Fig. \ref{fig:converge}, that, e.g., 
for the interval $x\in[-0.9,0.9]$ the linear convergence 
order is recovered.

\appendix
\section{Proof of Lemma \ref{lm:convergenza}}

%
%

\begin{proof}

We study convergence in two parts: the first concerns 
the accumulation of the error along the characteristic lines
of Eq. (\ref{basic_kolm_disc}); the second concern 
the error cumulated on the boundary integration 
of Eq. (\ref{basic_boundary_disc}).

Let $\ed$ be the global error (\ref{global_err})
between exact and discrete solutions.
By considering the standard procedure \cite{ran,mor:may} that puts
in relation local and global errors, we find from Eqs. (\ref{basic_kolm_charact}) 
and (\ref{basic_kolm_disc}):
\begin{equation}
\ed[k,j+1] = \ed -\Ad \alpha 
    \left( \ed[k+\nu,j] - \ed[k+\nu-1,j] \right) - \ld
    \ed \De y + \lu\varepsilon_{kj}^i \De y
\label{basic_kolm_err}
\end{equation} 
where $\lu\varepsilon_{kj}^i$ is defined in (\ref{loc_trunc}).
We start our analysis from Eq. (\ref{basic_kolm_err}), fixed
$i$ and given $j$ we find $\ed[k,j+1]$.
The error for $\ed[k,0]$ is found from the boundary integral 
condition.

The first step is to limit the error along the characteristic line $i$. 
From Eq. (\ref{basic_kolm_err}) we study the case 
$\Ad > 0$ ($\nu=0$):
$$
 |\ed[k,j+1]| \leq \left|1-\Ad \alpha - \ld \De y \right|
 | \ed | + \Ad \alpha |\ed[k-1,j]| + |\lu\varepsilon_{kj}^i| \De y
$$
Let $\Ed=\max_k|\ed|$ and $\lu\et_{j}^i=\max_k|\lu\varepsilon_{kj}^i|$ then:
$$
|\ed[k,j+1]|\leq \left(\left| 1-\Ad \alpha -\ld \De y  \right| 
                       + \Ad \alpha \right) \Ed + \lu\et_{j}^i \De y 
$$
If we set $1-\Ad \alpha - \ld \De y> 0 \;\;\;\; \forall l,j,k$,
i.e. the CFL condition, we get:
$$
|\ed[k,j+1]|\leq (1- \ld \De y ) \Ed + \lu\et_{j}^i \De y \;\;\;\;\;\; \forall k,j,l
$$
because of it is valid for all $k$, we can write:
\begin{equation}
	\Ed[j+1] \leq (1- \ld \De y) \Ed + \lu\et_{j}^i \De y
\label{major_E}
\end{equation} 
that is verified for $j=0,\ldots,n$ and $i=0,\ldots,n$.
Now we have to find $\lu E^{n+1}_0$.

The case $\nu=1$ gives the same result:
$$
 |\ed[k,j+1]| \leq |1+\Ad \alpha - \ld \De y|
 | \ed | - \Ad \alpha |\ed[k+1,j]|+|\lu\et_{j}^i| \De y
$$
and for the maximum error on $x$:
$$
|\ed[k,j+1]|\leq \left(\left| 1-|\Ad| \alpha -\ld \De y  \right| 
                       + |\Ad| \alpha \right) \Ed+ \lu\et_{j}^i \De y
$$
Let $1-|\Ad| \alpha - \ld \De y > 0$, we get Eq. (\ref{major_E}), so that 
it is valid independently by the sign of $A_l(x)$.

By iterating this expression we find:
\begin{equation}
    \Ed[m+1]\leq \Pi_{j=0}^m (1-\ld \De y) \Ed[0]+(m+1)\De y \lu\et^i
\label{bound_err_charact}
\end{equation} 
where $\lu\et^i = \max_j \lu\et_{j}^i$.
When the norm over all states is considered:
\begin{equation}
     \norm{E^i_{m}} \leq \norm{E^i_{0}} (1-L_u\De y)^m + m \De y
     \norm{\et^i}
\label{bound_err1_charat}
\end{equation}

Being $\ld \geq 0$ we have good chances to get an upper bound
to the error for increasing time.


Now we study the second step, i.e. the error $\ed[k,0]$ (or $\Ed[0]$)
along the boundary condition. Here $j=0$ then $n=i$.

From Eq. (\ref{basic_boundary_charact}) as calculated
with  quadrature on the exact solution $\phi_l(x,y,\xi)$, we can write:
$$
\phi_s(x_k,0,t_i) = 
  \sum_{l=1}^{S} q_{sl} \left( \sum_{j=0}^{i} \De y \w
  \lu \lambda(y_j) \phi_l(x_k,y_j,t_i-y_j) + \lu R_k^i \right)
$$
Subtracting side by side this equation and (\ref{basic_boundary_disc}):
$$
\lu[s]e_{k,0}^i = \sum_{l=1}^{S} q_{sl}\left(\sum_{j=0}^{i}\De y \w \ld
\lu e_{kj}^{i-j} + \lu R_k^i \right)
$$
and by applying the triangle inequality:
$$
|\lu[s]e_{k,0}^i|\leq\sum_{l=1}^{S} q_{sl}\left(\sum_{j=1}^{i}\De y \w \ld
|\lu e_{kj}^{i-j}| + |\lu R_k^i| \right).
$$
Let $\lu R^i = \max_k |\lu R_k^i|$ and being $q_{sl}\geq 0$,
$\w \geq 0$, we get:
\newcommand{\ldd}{\lu \underline{\lambda}}
\newcommand{\sEd}{\lu[s] E_0^i\;}
\begin{equation}
	\sEd \leq \sum_{l=1}^{S} q_{sl}\left(\sum_{j=0}^{i}\De y \w \ld
	\lu E_{j}^{i-j} + \lu R^i \right)
\label{bound_err_bound}
\end{equation} 

Let $\ldd = \min_j\; \ld$, from Eq. (\ref{bound_err_charact}) is:
$$
\lu E_{j+1}^i \leq (1-\ldd \De y)^{j+1} \Ed[0] + (j+1)\De y\lu\et^i
$$
and inserting it in (\ref{bound_err_bound}), we find:
\begin{equation}
\sEd \leq \sum_{l=1}^{S} q_{sl} \left(
  \De y \sum_{j=0}^{i} \w \ld [(1-\ldd \De y)^{j} \lu E_0^{i-j}+j\De y
  \lu\et^{i-j}] + \lu R^i
\right).
\end{equation} 
This inequality puts in relation the error on the boundary 
condition ($j=0$) with those at early times.

Let $\wu = \max_{j\geq 1} \w$ and $\ldu = \max_{j\geq 1}\; \ld$, then:
\begin{eqnarray}
   \sEd \leq \sum_{l=1}^{S} q_{sl}
     \left(
	\De y \wu \ldu \sum_{j=1}^i [(1-\ldd\De y)^j \lu E_0^{i-j} +
\right.&& j\De y   \lu\et^{i-j}]   \nonumber \\
&&     \left. +\De y w^{(i)}_0 \lu \lambda_0 \lu E_0^i +\lu R^i
     \right)   \nonumber
\end{eqnarray}
valid for $i\geq1$.
We have introduced the sum starting from $j=0$, so that, in order to find $\sEd$,
we have to solve a system of equation. 
However for convergence we do not need to solve it at all.
By moving $\lu E_0^{i}$ to the left hand side and summing over the states, we 
find:
\begin{eqnarray*}
    \sum_s \sEd- && \sum_{s,l} q_{sl}\De y w^{(i)}_0 \lu \lambda_0 \lu E_0^i 
     \leq  \\
    && \sum_s\sum_{l=1}^{S} q_{sl}
	\left(
	\De y \wu \ldu \sum_{j=1}^i [(1-\ldd\De y)^j \lu E_0^{i-j}+
	 j\De y   \lu\et^{i-j}]
 	+ \lu R^i
	\right)
\end{eqnarray*}
By using the fundamental property of stochastic matrix:
\newcommand{\Mwzero}{\bar{w}_0}
\begin{eqnarray*}
\sum_l (1- \Mwzero L_0\De y) && \lu E_0^i \leq  \\
&& \sum_{j=1}^i \De y \wu \sum_l \ldu [(1-\ldd\De y)^j \lu E_0^{i-j}
+ j\De y   \lu\et^{i-j}]+
\sum_l \lu R^i 
\end{eqnarray*}
where $\Mwzero=\max_i w_0^{(i)}$, $L_0=\max_l\lambda_0$. 
Let $a := 1-\Mwzero L_0 \De y >0$, we get:
$$
  \sum_l \lu E_0^i \leq  a^{-1} \De y
  \sum_{j=1}^i \wu \sum_l \ldu [(1-\ldd\De y)^j \lu E_0^{i-j}+
  j\De y   \lu\et^{i-j}] +\sum_l \lu R^i 
$$ 
where $0<a<1$.
Let $L_u = \max_l \ldu$, $L_d =\min_l \ldd$, $b=(1-L_d \De y)$ provided that from CFL
it is surely $L_d \De y \leq 1$, we have:
\begin{equation}
\norm{E_0^i} \leq a^{-1} \De y \wu L_u
\sum_{j=1}^i [b^j \norm{E_0^{i-j}} + j\De y\norm{\et^{i-j}}]+\norm{R^i}
\end{equation}
and with the summation index change $m=i-j$, it becomes:
$$
b^{-i} \norm{E_0^i} \leq \Ki \De y 
  \left(
     \De y \frac{i(i+1)}{2}\norm{\bar{\et}} +
     \sum_{m=0}^{i-1} b^{-m} \norm{E_0^{m}} 
  \right) + \\
  \norm{R^i} b^{-i}
$$
where $\norm{\bar{\et}}=\max_i \norm{\et^i}$, and $\Ki = \max_i(a^{-1}\wu L_u)$.
From the discrete Gronwall lemma (cf. in Th. 1.5.4 e Corol. 1.5.1 of \cite{bru}),
we obtain:
$$
b^{-i} \norm{E_0^i} \leq
\left( 
   \Ki \De y^2 \frac{i(i+1)}{2}\norm{\bar{\et}}+
   \norm{R^i}b^{-i}
\right)  (1+\Ki \De y)^i.
$$ 
Let $\norm{\bar{R}}=\max_i \norm{R^i}$, we find:
$$
\norm{E^i_0} \leq \norm{\bar{R}} (1+\Ki\De y)^i +
\norm{\bar{\et}} \Ki \De y^2 \frac{(i+1)i}{2} 
[b (1+\Ki \De y)]^i
$$
and  finally:
\begin{equation}
	\norm{E^i_0} \leq \norm{\bar{R}} \exp\left(\Ki t_i \right)   
 	+ \norm{\bar{\et}} \Ki \frac{t_{i+1}^2}{2}
	\exp\left(( \Ki - L_d)\, t_i \right)
\end{equation}

\end{proof}



\begin{thebibliography}{10}


\bibitem{bru}
{\sc H. Brunner and P.J. van der Houwen},
{\em The numerical solution of Volterra equation}, North-Holland
1986.

\bibitem{dav_bk}
{\sc M.H.A. Davis},
{\em Markov Models and Optimizations}, Chapman \& Hall/CRC, London 1993.

\bibitem{gar} 
{\sc C.~W. Gardiner},
{\em Handbook of Stochastic Methods}, Springer-Verlag Berlin Heidelberg New York, 
2nd Ed.~1985.

\bibitem{mor:may} 
{\sc K.~W. Morton and D.~F. Mayers},
{\em Numerical Solution of Partial Differential Equations}, 
Cambridge University Press 1994.

\bibitem{ran} 
{\sc Randall J.~LeVeque}, 
{\em Numerical Methods for conservation Laws},
Birkh\"auser 1992.


\bibitem{zwa}
{\sc R.~Zwanzig},
{\em Lectures in Theoretical Physics (Boulder)},
vol. 3 ed. W.E.~Broton (New York: Interscience) 1961.


\bibitem{mario:prep}
{\sc M.~Annunziato},
{\em A finite difference method for piecewise deterministic Markov processes},
e-print: {\tt http://arxiv.org/abs/math.NA/0606588}.


\bibitem{mario} 
{\sc M.~Annunziato}, 
{\em Non-gaussian equilibrium distribution arising from the Langevin equation}, 
Phys. Rev. E, {65} (2002), pp.~21113-1--6.



\bibitem{grigo}
{\sc M.~Annunziato, P.~Grigolini and J.~Riccardi},
{\em Fluctuation-dissapation process without a time scale},
Phys. Rev. E, {61} (2000), pp.~4801--4808.
{\sc M.~Annunziato, P.~Grigolini and B.~J. West},
{\em Canonical and noncanonical equilibrium distribution},
Phys. Rev. E, {64} (2001), pp.~011107-1--13.


\bibitem{ben:etal}
{\sc I.~Bena, C.~Van den Broeck, R.~Kawai and K.~Lindenberg},
{\em Drift by dichotomous Markov noise},
Phys. Rev. E, {68} (2003), pp.~041111-1--12.

\bibitem{bol:gri:wes} 
{\sc M.~Bologna, P.~Grigolini, and B.~J. West}, 
{\em Strange kinetics: conflict between density and trajectory description}, 
Chem. Phys., {284} (2002), pp.~115--128.

\bibitem{bud:cac} 
{\sc A.~Budini and M.~O C\'aceres}, 
{\em Functional characterization of generalized Langevin equations},
J. Phys. A: Math. Gen., {37} (2004), pp.~5959--5981.

\bibitem{cac}
{\sc M.~O C\'aceres}, 
{\em Computing a non-Maxwellian velocity distribution from first principles},
Phys. Rev E, {67} (2003), pp.~016102--016104.
 

\bibitem{cos:duf}
{\sc O.~L.V.~Costa  and F.~Dufour},
{\em On the Poisson equation for piecewise-deterministic Markov processes},
SIAM J.Control Optim., {42} (2003), pp.~985--1001.

\bibitem{cox2}
{\sc D.R.~Cox},
{\em The analysis of non-Markovian stochastic processes by the inclusion
of supplementary variables},
Proc. Camb. Phil. Soc., {51} (1955), pp.~433--441.


\bibitem{dav}
{\sc M.~H.A. Davis},
{\em Piecewise-Deterministic Markov Processes: A General Class of Non-Diffusion Stochastic Models}, J. of the Royal Stat. Soc. Series B, {46} (1984), pp.~353--388.

\bibitem{fil:hon}
{\sc R.~Filliger and M.~O. Hongler}, 
{\em Supersymmetry in random two-velocity processes},
 Physica A, {332} (2004), pp.~141--150.

\bibitem{hor} 
{\sc W.~Horsthemke}, 
{\em Spatial instabilities in reaction random walks with direction-independent kinetics}, 
Phys. Rev E, {60} (1999), pp.~2651--2663.


\bibitem{jak:ren} 
{\sc E.~Jakeman and E.~Renshaw}, 
{\em Correlated random-walk model for scattering},
J. Opt. Soc. Am. A, {4} (1987), pp.~1206--1212.


\bibitem{jak:rid} 
{\sc E.~Jakeman and  K.~D. Ridley},
{\em The statistics of a filtered telegraph signal},
J. Phys. A, {32}  (1999), pp.~8803--8821.

\bibitem{kit} 
{\sc K.~Kitahara \textit{et al.}}, 
{\em Phase Diagrams of Noise Induced Transitions},
Prog. Theor. Phys., {64} (1980), pp.~1233--1247.


\bibitem{man:etal}
{\sc R.~Mankin, A.~Ainsaar, A.~Haljas and E.~Reiter},
{\em Trichotomous-noise-induced catastrophic shifts in symbiotic ecosystems},
Phys. Rev. E, {65} (2002), pp.~051108-1--9.



\bibitem{mas:por} 
{\sc J.~Masoliver and J.~M. Porr\`a}, 
{\em Harmonic oscillators driven by colored noise: Crossovers, resonances, and spectra}, 
Phys. Rev. E, {48} (1993), pp.~4309--4319.



\bibitem{mori.h}
{\sc H.~Mori},
{\em Statistical-Mechanical Theory of Transport in Fluids},
Phys. Rev., {112} (1958), pp.~1829--1842.



\bibitem{mor} 
{\sc A.~Morita},
{\em Free Brownian motion of a particle driven by a dichotomous random force}, 
Phys. Rev. A, {41} (1990), pp.~754--760.



\bibitem{hil:oth} 
{\sc T.~Hillen and H.~G. Othmer}, 
{\em The diffusion limit of transport equations derived from velocity-jump processes},
SIAM J. Appl. Math., {61} (2000), pp.~751--775.


\bibitem{oth:dun:alt} 
{\sc H.~G. Othmer, S.~R. Dunbar, and W.~Alt}, 
{\em Models of dispersal in biological systems},
J. Math. Biol., {26} (1988), pp.~263--298.



\bibitem{paw:ric} 
{\sc R.~F. Pawula and O.~Rice},
{\em On Filtered Binary Processes}, 
IEEE Trans. Inf. Th., Vol. {IT-32} (1986), pp.~63--72.

\bibitem{paw}
{\sc R.~F. Pawula}, 
{\em Approximating distributions from moments},
Phys. Rev. A, {36} (1987), pp.~4996--5007.


\bibitem{hyper}
{\sc W.~F. Perger, A.~Bhalla and M.~Nardin}, 
{\em A numerical evaluator for the generalized hypergeometric 1 series},
Comp. Phys. Commun., {77} (1993), pp.~249--254
(MATLAB\circledR \, translation by B.~E. Barrowes, see MATLAB\circledR \, Central File Exchange).
This routine uses a gamma function for complex arguments, see {\sc P.~Godfrey},
{\tt http://my.fit.edu/$\sim$gabdo/gamma.m}.



\bibitem{bro:han} 
{\sc C.~Van den Broeck and P.~H\"anggi}, 
{\em Activation rates for nonlinear stochastic flows driven by non-Gaussian noise}, Phys. Rev. A, {30} (1984), pp.~2730--2736.
 

\bibitem{won} 
{\sc W.~M. Wonham},  
J. Electron. Control, {6} (1959), pp.~376--383. 




\end{thebibliography}
\end{document}